\documentclass[12pt]{amsart}

\usepackage{amssymb}
\usepackage{latexsym}
\def\C{\mathbb C}

\def\N{\mathbb N}
\def\R{\mathbb R}
\def\FF{{\mathcal F}}

\def\re{\operatorname{Re}}

\newtheorem{la}{Lemma}
\newtheorem*{tha}{Theorem A}
\newtheorem*{thb}{Theorem B}

\newtheorem{thm}{Theorem}

\theoremstyle{definition}
\newtheorem*{question}{Problem}
\newtheorem*{conja}{Conjecture A}
\newtheorem*{conjb}{Conjecture B}
\newtheorem*{conjc}{Conjecture C}

\begin{document}
\title{Normal families and fixed points of iterates}
\dedicatory{Dedicated to Professor Yang Lo on the occasion
of his 70th birthday}
\author{Walter Bergweiler}
\thanks{Supported by the EU Research Training Network CODY, the 
ESF Networking Programme HCAA 
and the Deutsche Forschungsgemeinschaft, Be 1508/7-1.}
\address{Mathematisches Seminar der 
Christian--Albrechts--Universit\"at zu Kiel,
Ludewig--Meyn--Stra{\ss}e~4,
D--24098 Kiel,
Germany}
\email{bergweiler@math.uni-kiel.de}
\date{\today}
\begin{abstract}
Let $\FF$ be a family of holomorphic functions and suppose
that there exists
$\varepsilon>0$ such that if $f\in \FF$, then
$|(f^2)'(\xi)|\leq 4-\varepsilon$ for all fixed points $\xi$ 
of the second iterate~$f^2$. We show that then $\FF$ is normal.
This is deduced from a result which says that if
$p$ is a polynomial of degree at least $2$, 
then $p^2$ has a fixed point 
$\xi$ such that $|(p^2)'(\xi)|\geq 4$. 
The results are motivated by a problem posed by Yang Lo.
\end{abstract}
\maketitle

\section{Introduction and main results}
Yang Lo~\cite[Problem~8]{Yan92} posed the following problem in 1992.
\begin{question}
Let $\FF$ be a family of entire functions, let
$D\subset \C$ be a domain and let $n\geq 2$ be a fixed integer.
Suppose that for every $f\in\FF$ the $n$-th iterate $f^n$
does not have fixed points in $D$.  Is $\FF$ normal in~$D$?
\end{question}
An affirmative answer was given by  Ess\'en and Wu~\cite{Ess98}
in 1998. 
They did not require that the functions in $\FF$ are 
entire but only that they are holomorphic in~$D$.
The iterates $f^n:D_n\to \C$ of such a function $f:D\to \C$ 
are defined by $D_1:=D$, $f^1:=f$ and
$D_n:=f^{-1}(D_{n-1})$, 
$f^n:=f^{n-1}\circ f$ for $n\in\N$, $n\geq 2$.
Note that $D_2=f^{-1}(D_1)\subset D=D_1$ and thus $D_{n+1}\subset
D_n\subset D$ for all $n\in\N$.

In a subsequent paper, Ess\'en and Wu~\cite[Theorem~1]{Ess00} 
gave the following generalization of their result.
Here a fixed point $\xi$ of $f$ is called {\em repelling} 
if $|f'(\xi)|>1$. 
\begin{tha}
Let $D\subset \C$ be a domain and let $\FF$ be the family of all
holomorphic functions $f:D\to\C$ for which there exists $n=n(f)>1$ 
such that $f^n$ has no repelling fixed point. Then $\FF$ is normal.
\end{tha}
There are a number of further developments initiated by 
Yang Lo's question.
For example, 
his question has also been considered 
for meromorphic~\cite{Wan02} and
quasiregular~\cite{Sie06} functions.
Other papers related to Yang Lo's problem include~\cite{BB,Cha05,Cha08,Xu05};
see~\cite[section 3]{Ber06} for further discussion.
The following result was proved in~\cite[Theorem~1.3]{Ber}.
\begin{thb} 
For each integer 
$n\geq 2$
there exists a constant
$K_n>1$ with the following property:
if $D\subset \C$ is a domain and $\FF$ is a family of
holomorphic functions $f:D\to\C$ such that $|(f^n)'(\xi)|\leq K_n$
for all fixed points $\xi$ of $f^n$,
then $\FF$ is normal.
\end{thb}
Considering the family $\FF=\{az^2\}_{a\in\C\backslash\{0\}}$
we see that the conclusion of Theorem~B does not hold 
for $K_n=2^n$. The following conjecture says that it holds 
for $K_n<2^n$.
\begin{conja} 
Let $D\subset\C$ be a domain, $n\geq 2$ and $\varepsilon>0$.
Let $\FF$ be the family of all
holomorphic functions $f:D\to\C$ such that 
$|(f^n)'(\xi)|\leq 2^n-\varepsilon$
for all fixed points $\xi$ of $f^n$.
Then $\FF$ is normal.
\end{conja} 
We show that this conjecture is true for $n=2$.  
\begin{thm} \label{thm1} 
Let $D\subset\C$ be a domain and $\varepsilon>0$.
Let $\FF$ be the 
family of all holomorphic functions $f:D\to\C$ such 
that $|(f^2)'(\xi)|\leq 4-\varepsilon$  for all fixed 
points $\xi$ of $f^2$.  Then $\FF$ is normal.
\end{thm}
We deduce Theorem~\ref{thm1} from a result about fixed points
of iterated polynomials. 
In fact, we shall see that Conjecture A is equivalent to the following
conjecture.
\begin{conjb} 
Let $p$ be a polynomial of degree at least $2$ and let  $n\geq 2$.
Then $p^n$ has a fixed point $\xi$ satisfying $|(p^n)'(\xi)|\geq 2^n$.
\end{conjb} 
The equivalence of these two conjectures is seen by the following 
result.
\begin{thm} \label{thm2}
Let $n\geq 2$ and 
let $C_n>0$ be such that for every  polynomial $p$ of degree at least $2$
there exists a fixed point $\xi$ of $p^n$ such that $|(p^n)'(\xi)|\geq C_n$.
Let $D\subset \C$ be a domain, let 
$\varepsilon>0$
 and let $\FF$ be the family of all
holomorphic functions $f:D\to\C$ such that $|(f^n)'(\xi)|\leq C_n-
\varepsilon$ for all fixed points $\xi$ of $f^n$.
Then $\FF$ is normal.
\end{thm}
Theorem~\ref{thm1} now follows from Theorem~\ref{thm2} and the
following result which says that we can take $C_2=4$.
\begin{thm} \label{thm3}
Let $p$ be a polynomial of degree at least $2$.
Then $p^2$ has a fixed point $\xi$ satisfying $|(p^2)'(\xi)|\geq 4$.
\end{thm}
We conclude this introduction with a conjecture which is stronger than
Conjecture~B.
\begin{conjc} 
Let $p$ be a polynomial of degree $d\geq 2$ and let  $n\geq 2$.
Then $p^n$ has a fixed point $\xi$ satisfying $|(p^n)'(\xi)|\geq d^n$.
\end{conjc} 
The monomial $p(z)=z^d$ shows that this would be best possible.
A.~E.~Eremenko and G.~M.~Levin~\cite[Theorem~3]{EreLev}
have shown that if $p$ is a polynomial of degree $d\geq 2$ which is
not conjugate to the monomial $z^d$, then there exists $n\geq 2$ such
that  $p^n$ has a fixed point $\xi$ satisfying $|(p^n)'(\xi)|> d^n$.

\section{Proof of Theorem~\ref{thm2}}
We shall use the following result proved in~\cite[Theorem~1.2]{Ber}.
\begin{la}\label{trans}
Let $f$ be a transcendental entire function and let $n\in\N$, $n\geq 2$.
Then there exists a sequence  $(\xi_k)$ of fixed points of $f^n$
such that $(f^n)'(\xi_k)\to\infty$ as $k\to\infty$.
\end{la}
The other main tool in the proof of Theorem \ref{thm2} is the following lemma
due to X.~Pang and L.~Zalcman~\cite[Lemma 2]{PZ}.
\begin{la} \label{PangZalcman}
Let $\FF$ be a family of functions meromorphic on the unit disc,
all of whose zeros have multiplicity at least $k$, and suppose
that there exists $A\geq 1$ such that $|g^{(k)}(\xi)|\leq A$
whenever $g(\xi)=0$, $g\in\FF$. Then if $\FF$ is not normal there
exist, for each $0\leq \alpha\leq k$, a number $r\in (0,1)$,
points $z_j\in D(0,r)$, functions $g_j\in \FF$ and positive
numbers $\rho_j$ tending to zero such that
$$\frac{g_j(z_j+\rho_jz)}{\rho_j^\alpha}\to G(z)$$
locally uniformly, where $G$ is a nonconstant
meromorphic function on $\C$
such that the spherical derivative
$G^\#$ of $G$ satisfies
$G^\#(z)\leq G^\#(0)=kA+1$ for all $z\in \C$.
\end{la}
We shall only need the case $k=1$ of Lemma~\ref{PangZalcman}.
This special case can also be found in Pang's paper~\cite[Lemma~2]{Pang02}.
The case $\alpha=0$ is known as Zalcman's lemma~\cite{Zal75,Zal98}.

\begin{proof}[Proof of Theorem \ref{thm2}]
We denote by $\FF(D,n,K)$ the family 
of all functions $f$ holomorphic in $D$ such
that $|(f^n)'(\xi)|\leq K$ whenever  $f^n(\xi)=\xi$.
Note that this implies that 
$|f'(\xi)|\leq \sqrt[n]{K}$ whenever  $f(\xi)=\xi$.
Suppose that the conclusion of Theorem~\ref{thm2} does not hold.
Then there exist a domain $D\subset \C$, 
$n\geq 2$ and $\varepsilon>0$ 
such that $\FF(D,n,C_n-\varepsilon)$ is not normal.
We may assume that $D$ is the unit disk.

We choose
a non-normal sequence $(f_j)$ in
$\FF(D,n,C_n-\varepsilon)$. With $g_j(z):=f_j(z)-z$ we find that
if $g_j(\xi)=0$, then $f_j(\xi)=\xi$ and thus
$$|g_j'(\xi)|\leq |f_j'(\xi)|+1\leq \sqrt[n]{C_n-\varepsilon}+1=:A.$$
We may assume here that $\varepsilon$ is chosen so small that $A^n>C_n$.
Clearly, the sequence $(g_j)$ is also not normal.
Applying Lemma~\ref{PangZalcman} with $\alpha=k=1$
we may assume,
passing to a subsequence if necessary, that there exist $z_j\in D$
and $\rho_j>0$ such that $g_j(z_j+\rho_j z)/\rho_j\to G(z)$ for
some entire function $G$ satisfying 
$G^\#(z)\leq  G^\#(0)=A+1$ 
for all $z\in\C$.
With $L_j(z)=z_j+\rho_j z$ we find that
\[
h_j(z):=L_j^{-1}(f_j(L_j(z)))=\frac{f_j(z_j+\rho_j z)-z_j}{\rho_j}
=\frac{g_j(z_j+\rho_j z)}{\rho_j}+z\to G(z)+z.
\]
With $F(z):=G(z)+z$ we thus have $h_j(z)\to F(z)$ 
as $j\to\infty$. It follows that
$h_j^n(z)
\to F^n(z)$.
The assumption that
$f_j\in \FF(D,n,C_n-\varepsilon)$ implies that 
$h_j\in \FF(L_j^{-1}(D),n,C_n-\varepsilon)$; that is,
$|(h_j^n)'(\xi)|\leq C_n-\varepsilon$ whenever
$h_j^n(\xi)=\xi$. We deduce that $F\in \FF(\C,n,C_n-\varepsilon)$.
It follows from the definition of $C_n$ that $F$ cannot be a 
polynomial of degree greater than one. 
And Lemma~\ref{trans} implies that $F$ cannot be transcendental.
Thus $F$ is a polynomial of degree $1$ at most.
Now $|F'(0)|\geq |G'(0)|-1\geq G^\#(0)-1=A.$
Hence $F$ has the form
$F(z)=az+b$ where $|a|\geq A$.
With $\xi:=b/(1-a)$ we obtain $F(\xi)=\xi$ and
$|F'(\xi)|=|a|\geq A$.
Thus $F^n(\xi)=\xi$ and
$|(F^n)'(\xi)|=|a^n|\geq A^n>C_n$,
contradicting $F\in \FF(\C,n,C_n-\varepsilon)$.
\end{proof}
\section{Proof of Theorem~\ref{thm3}}
The following lemma is due to A.~E.~Eremenko and 
G.~M.~Levin~\cite[Lemma~1]{EreLev}.
\begin{la}\label{EL}
Let $p$ be a polynomial of degree $d\geq 2$.
Then there exists $c$ such that 
$$\sum_{\{z:p^n(z)=z\}} (p^n)'(z)=d^n(d^n-1)+c^n$$
for all $n\in\N$.
\end{la}
In fact, they show that this holds with 
$$c=\sum_{\{z:p(z)=w\}} p'(z),$$
for arbitrary $w\in\C$, but we do not need this result.
We note that in the sum occuring in the lemma
each fixed point of $p^n$ is counted 
according to multiplicity.

\begin{proof}[Proof of Theorem \ref{thm3}]
Suppose that $p$ is a polynomial such that $|(p^2)'(\xi)|< 4$
for each fixed point $\xi$ of $p^2$. Then
$|p'(\xi)|< 2$
for each fixed point $\xi$ of $p$.
It follows from Lemma~\ref{EL}
that 
\begin{equation}\label{el1}
|d(d-1)+c|<2d
\end{equation}
and 
\begin{equation}\label{el2}
|d^2(d^2-1)+c^2|<4d^2.
\end{equation}
Now (\ref{el1}) yields
$c=-d(d-1)+re^{it}$ where $0\leq r<2d$ and $t\in\R$.
Thus
$$c^2=(-d(d-1)+re^{it})^2=d^2(d-1)^2-2d(d-1)re^{it}+r^2e^{i2t}$$
and hence
\begin{eqnarray*}
\re(c^2)
&=&
d^2(d-1)^2-2d(d-1)r\cos t +r^2\cos 2t\\
&=&
d^2(d-1)^2-2d(d-1)r\cos t +r^2(2\cos^2 t -1)\\
&=&
d^2(d-1)^2-r^2 -2d(d-1)r\cos t + 2r^2 \cos^2 t\\
&=&
\frac12 d^2(d-1)^2 -r^2+
2\left(\frac12 d(d-1)-r\cos t \right)^2 \\
&\geq&
\frac12 d^2(d-1)^2-r^2\\
&>&
\frac12 d^2(d-1)^2-4d^2.
\end{eqnarray*}
Thus 
\begin{eqnarray*}
|d^2(d^2-1)+c^2|&\geq& \re(d^2(d^2-1)+c^2)\\
&>&
d^2(d^2-1)+\frac12 d^2(d-1)^2-4d^2\\
&=&d^2\left(\frac32 d^2 -d-\frac92\right).
\end{eqnarray*}
Combining this with (\ref{el2}) we find that
$$d^2\left(\frac32 d^2 -d-\frac92\right)<4d^2$$
and thus that
$$\frac32 d^2 -d-\frac92<4.$$
The last inequality can be rewritten as
$$\frac32 d^2 -d -\frac{17}{2}=\frac32\left( 
\left(d-\frac13\right)^2- \frac{52}{9}
\right)<0.$$
It follows that 
$d<(1+2\sqrt{13})/3<3$.

It thus remains to consider the case $d=2$. 
Then $c=0$ in Lemma~\ref{EL}. Let 
$\xi_1,\xi_2$ be the fixed points of $p$ and
$\xi_1,\xi_2,\xi_3,\xi_4$ those of  $p^2$.
Then $p'(\xi_1)+p'(\xi_2)=2$ so that
$p'(\xi_{1,2})=1\pm a$ for some $a\in\C$,
with $|1\pm a|<2$.
In particular, we have $|\re a|<1$.
Moreover, Lemma~\ref{EL} yields
$$12=\sum_{j=1}^4 (p^2)'(\xi_j)=(1+a)^2+(1-a)^2+(p^2)'(\xi_3)+(p^2)'(\xi_4)$$
and hence
$$(p^2)'(\xi_3)+(p^2)'(\xi_4)=2(5-a^2).$$
It follows that 
$$\re\left( (p^2)'(\xi_3)+(p^2)'(\xi_4)\right)=2(5-\re( a^2))
\geq 2(5-(\re a)^2)>8.$$
Hence there exists $j\in\{3,4\}$ with
$|(p^2)'(\xi_j)|>4$, a contradiction.
\end{proof}

\end{document}